\documentclass[12pt]{amsart}
\usepackage[colorlinks=true,
pdfstartview=FitV, linkcolor=blue, citecolor=red,
urlcolor=blue]{hyperref}
\usepackage{amscd,amssymb,epsf}
\usepackage[mathscr]{eucal}
\newtheorem{thm}{Theorem}
\newtheorem*{thm*}{Theorem}
\newtheorem{lemma}[thm]{Lemma}

\newtheorem*{conjecture*}{Conjecture}
\newtheorem*{question*}{Question}
\newtheorem{proposition}[thm]{Proposition}
\newtheorem{corollary}[thm]{Corollary}
\newtheorem*{corollary*}{Corollary}
\newtheorem*{claim*}{Claim}

\numberwithin{equation}{section}
\numberwithin{thm}{section}
\newtheorem*{OSLemma*}{Opposite Sign Lemma}

\begin{document}

\title[Proper affine actions]
{Proper affine actions and geodesic flows of hyperbolic surfaces}
\author[Goldman]{W.  M. Goldman}
\address{ Department of Mathematics,
University of Maryland, College Park, MD  20742 USA  }
\email{ wmg@math.umd.edu }
\author[Labourie]{F. Labourie}
\address{Topologie et Dynamique, 
Universit\'e Paris-Sud F-91405 Orsay (Cedex)
France}
\email{francois.labourie@math.u-psud.fr}
\author[Margulis]{G. Margulis}
\address{Department of Mathematics,
10 Hillhouse Ave.,  
P.O. Box 208283, 
Yale University,
New Haven, CT 06520 USA}
\email{margulis@math.yale.edu}

\thanks{Goldman 
gratefully acknowledge partial support from 
National Science Foundation grants  DMS-0103889,
DMS0405605, DMS070781, the Mathematical Sciences
Research Institute and the Oswald Veblen Fund
at the Insitute for Advanced Study, and a Semester
Research Award from the General Research Board
of the University of Maryland.
Margulis gratefully acknowledge partial support from 
National Science Foundation grant  DMS-0244406.
Labourie thanks l'Institut Universitaire de France.}

\date{\today}

\newcommand\R{\mathbb{R}}
\newcommand\Z{\mathbb{Z}}
\newcommand\B{\mathbb{B}}
\newcommand\bB{\mathsf{B}}
\renewcommand\L{\mathbb{L}}
\newcommand\E{\mathbb{E}}
\newcommand\eE{\mathsf{E}}
\newcommand\V{\mathbb{V}}
\newcommand\vV{\mathsf{V}}
\newcommand\G{\mathsf{G}}
\newcommand\Go{\mathsf{G}_0}
\renewcommand\gg{\mathsf{g}}
\newcommand\Rto{\R^{2,1}}
\newcommand\Eto{\E^{2,1}}
\newcommand\Ht{\mathsf{H}^2}
\newcommand\Id{\mathbb{I}}
\newcommand\Isom{\mathsf{Isom}}
\newcommand\Isomo{\Isom^0(\Eto)}
\newcommand\Oto{\mathsf{O}(2,1)^0}
\newcommand\Ot{\mathsf{O}(2,1)}
\newcommand\Aff{\mathsf{Aff}}
\newcommand\A{\mathsf{A}_0}
\newcommand\GL{\mathsf{GL}}
\newcommand\SL{\mathsf{SL}(2,\R)}
\newcommand\SO{\mathsf{SO}(2)}

\newcommand\GLV{\GL(\vV)}
\newcommand\AffE{\Aff(\eE)}
\newcommand\GLn{\GL(\R^n)}

\newcommand\GLthr{\GL(\R^3)}
\newcommand\Ee{\mathcal{E}}
\renewcommand\Pr{\mathcal{P}}
\newcommand\kr{\mathcal{K}_R}
\newcommand\pp{\mathsf{p}_0}  %%%% maybe change?
\newcommand\qq{\mathsf{q}}
\newcommand\ppn{\mathsf{p}_n}
\newcommand\qqn{\mathsf{p}_n}
\newcommand\uu{\mathsf{u}}
\newcommand\vv{\mathsf{v}}
\newcommand\ww{\mathsf{w}}
\renewcommand\ss{\mathsf{s}}
\newcommand\tvv{\tilde{\mathsf{v}}}
\newcommand\F{F_{\uu,s}}
\newcommand\PSL{\mathsf{PSL}(2,\R)}
\newcommand\homeo{{\mathsf{Homeo}(X)}}
\newcommand\PGL{\mathsf{PGL}(2,\R)}

\newcommand\so{S^1_\infty}
\newcommand\Ur{U_{\mathsf{rec}}\Sigma}
\newcommand\tUr{\widetilde{\Ur}}
\newcommand\Er{\E_{\mathsf{rec}}\Sigma}
\newcommand\Vr{\V_{\mathsf{rec}}\Sigma}
\newcommand\Vrp{\V^{+}_{\mathsf{rec}}\Sigma}
\newcommand\Vrm{\V^{-}_{\mathsf{rec}}\Sigma}
\newcommand\Vrpm{\V^{\pm}_{\mathsf{rec}}\Sigma}
\newcommand\ho{H^1(\Gamma_0,\Rto)}

\newcommand\K{\mathsf{K}_0}
\newcommand\Ss{\mathcal{S}(\E)}
\newcommand\gc{\mathcal{C}(\Sigma)}
\newcommand\Az{\mathcal{R}^0(\V)}
\newcommand\U{{U\Sigma}}
\newcommand\T{\mathbb T}
\newcommand\xo{\mathsf{x}^0}
\newcommand\Vm{\mathsf{V}^-}
\newcommand\Vp{\mathsf{V}^+}
\newcommand\Vpm{\mathsf{V}^\pm}
\newcommand\Vmp{\mathsf{V}^\mp}

\newcommand\tSs{\mathcal{S}(\tE)}
\newcommand\tAz{\mathcal{R}^0(\tV)}
\newcommand\tPi{\tilde{\Pi}}
\newcommand\tPhi{\tilde{\Phi}}

\newcommand\z{Z^1(\Gamma_0,\vV)}
\newcommand\h{H^1(\Gamma_0,\vV)}
\newcommand\Ph{\mathbb{P}\big(\h\big)}

\newcommand\gcp{\mathcal{C}_{\mathsf{per}}(\Sigma)}
\newcommand\grp{\mathcal{G}}
\newcommand\compact{\subset\hspace{-2pt}\subset}
\newcommand{\lt}{\Lambda^{\star}_2}
\renewcommand{\ker}{\mathsf{Ker}}
\newcommand{\tv}{\tilde{v}}
\newcommand{\sv}{\mathsf{v}}
\newcommand{\tnu}{\tilde{\nu}}
\newcommand{\ts}{\tilde{s}}
\newcommand{\tx}{\tilde{x}}
\newcommand{\tnabla}{\tilde{\nabla}}
\newcommand{\tphi}{\tilde{\phi}}
\newcommand{\Ddt}{\frac{D}{dt}}
\newcommand{\ddt}{\frac{d}{dt}}
\newcommand{\tV}{\tilde{\V}}
\newcommand{\tE}{\tilde{\E}}
\newcommand{\core}{{\mathsf{core}(\Sigma)}}
\newcommand{\supp}{\operatorname{supp}(\mu)}
\newcommand{\xmapsto}[1]{|\hspace{-5pt}\xrightarrow{~{#1}~}}

\renewcommand{\phi}{\varphi} 
%% only use varphi to avoid confusion with upper case \Phi
%% 

\setcounter{tocdepth}{3} %% change this to 3 for editing purposes

\begin{abstract}
Let $\Gamma_0\subset\Ot$ be a Schottky group, and let
$\Sigma=\Ht/\Gamma_0$ be the corresponding hyperbolic surface.  Let
$\gc$ denote the space of geodesic currents on $\Sigma$.
The cohomology group $\h$ parametrizes
equivalence classes of affine deformations 
$\Gamma_\uu$ of $\Gamma_0$ acting on an irreducible representation $\vV$ of 
$\Ot$. 
We define a continuous biaffine map
\begin{equation*}
\gc \times \h \xrightarrow{\Psi} \R 
\end{equation*}
which is linear on the vector space $\h$.
An affine deformation $\Gamma_\uu$ acts properly
if and only if 
$\Psi(\mu,[\uu])\neq 0$
for all $\mu\in\gc$. Consequently the set of proper affine actions 
whose linear part is a Schottky group identifies with a bundle
of open convex cones in $\h$ over the Teichm\"uller space of $\Sigma$.
\end{abstract}

\maketitle
\pagebreak
\tableofcontents
\pagebreak
\section*{Introduction}

It is well known that every discrete group of Euclidean isometries of
$\R^n$ contains a free abelian subgroup of finite index.
In~\cite{Milnor}, Milnor asked if every discrete subgroup of {\em
affine transformations \/} of $\R^n$ must contain a polycyclic
subgroup of finite index.  Furthermore he showed that this question is
equivalent to whether a discrete subgroup of $\Aff(\R^n)$ which acts
properly must be {\em amenable.\/} By Tits~\cite{Tits}, this is
equivalent to the existence of a proper affine action of a nonabelian
free group.  Margulis subsequently showed~\cite{Margulis1,Margulis2}
that proper affine actions of nonabelian free groups do indeed
exist. 
%%%
%%%%The present paper completes the classification of proper affine
%%%actions on $\R^3$ for a large class of discrete groups. 
%%%
The present paper describes the deformation space of proper affine
actions on $\R^3$ for a large class of discrete groups, as a bundle
of convex cones over the Fricke-Teichm\"uller deformation space of
hyperbolic structures on a compact surface with geodesic boundary.
%%%

If $\Gamma\subset\Aff(\R^n)$ is a discrete subgroup which acts properly on
$\R^n$, then a subgroup of finite index will act freely. In this case the
quotient $\R^n/\Gamma$ is a {\em complete affine manifold $M$\/} 
with fundamental group $\pi_1(M)\cong \Gamma$.
When $n=3$, Fried-Goldman~\cite{FriedGoldman} classified all
the amenable such $\Gamma$ (including the case when $M$ is compact).
Furthermore the linear holonomy homomorphism 
\begin{equation*}
\pi_1(M) \xrightarrow{\L} \GLthr 
\end{equation*}
embeds $\pi_1(M)$ onto a discrete subgroup 
$\Gamma_0\subset\GLthr$, which  is conjugate to a
subgroup of $\Oto$ (Fried-Goldman~\cite{FriedGoldman}). Mess~\cite{Mess} 
proved the hyperbolic surface
\begin{equation*}
\Sigma = \Gamma_0\backslash\Oto/\SO 
\end{equation*}
is noncompact.  
Goldman-Margulis~\cite{GoldmanMargulis} and Labourie~\cite{Labourie}
gave alternate proofs of Mess's theorem, using ideas which evolved to 
the present work. 
% Labourie's result applies to affine deformations
% of Fuchsian actions in all dimensions.  
% moved

Thus the classification of proper 3-dimensional affine actions reduces
to {\em affine deformations\/} of free discrete subgroups
$\Gamma_0\subset\Oto$.  An {\em affine deformation \/} of $\Gamma_0$
is a group $\Gamma$ of affine transformations whose linear part equals
$\Gamma_0$, that is, a subgroup $\Gamma\subset\Isomo$ such that the
restriction of $\L$ to $\Gamma$ is an isomorphism
$\Gamma\longrightarrow\Gamma_0$.

Equivalence classes of affine
deformations of $\Gamma_0\subset\Oto$
are para\-metrized by the cohomology group $\h$, where $\vV$ is the linear
holonomy representation. Given a cocycle $\uu\in\z$, we denote
the corresponding affine deformation by $\Gamma_{\uu}$.
Drumm~\cite{Drumm1,Drumm2,Drumm4,CharetteGoldman}
showed that Mess's necessary condition of non-compactness is sufficient: 
{\em every noncocompact discrete subgroup of $\Oto$ admits 
proper affine deformations.\/} In particular he found an open subset
of $\h$ parametrizing proper affine deformations~\cite{DrummGoldman}.

This paper gives a criterion for the properness of 
an affine deformation $\Gamma_\uu$ in terms of the parameter $[\uu]\in\h$.

With no extra work, we take $\vV$ to be a representation of $\Gamma_0$
obtained by composition of the discrete embedding $\Gamma_0\subset\SL$
with any irreducible representation of $\SL$.  Such an action is
called {\em Fuchsian\/} in Labourie~\cite{Labourie}.  Proper actions
occur only when $\vV$ has dimension $4k+3$ (see Abels~\cite{Abels},
Abels-Margulis-Soifer~\cite{AbelsMargulisSoifer,AbelsMargulisSoifer2}
and Labourie~\cite{Labourie}).  In those dimensions, 
Abels-Margulis-Soifer~\cite{AbelsMargulisSoifer,AbelsMargulisSoifer2}
constructed proper affine deformations of Fuchsian actions of free groups.

% One of the main results of the present paper is the following:

\begin{thm*} Suppose that $\Gamma_0$ contains no parabolic elements.
Then the equivalence classes of proper affine deformations of $\Gamma$
form an open convex cone in $\h$.
\end{thm*}

The main tool in this paper is a generalization of 
the invariant constructed by 
Margulis~\cite{Margulis1,Margulis2} and extended to higher dimensions
by Labourie~\cite{Labourie}.
%% moved from page 2
%% identify ``Margulis invariant''
Margulis's invariant $\alpha_\uu$ is
a class function
$\Gamma_0 \stackrel{\alpha_\uu}\longrightarrow \R$ 
associated to an affine deformation $\Gamma_\uu$.
It satisfies the following properties:
\begin{itemize}
\item
$\alpha_{[\uu]}(\gamma^n) = \vert n\vert \alpha_{[\uu]}(\gamma)$
%% (when $r$ is odd) %% 
;  
\item
$\alpha_{[\uu]}(\gamma) = 0$ $\Longleftrightarrow$ 
$\gamma$ fixes a point in 
$\V$; 
\item
The function $\alpha_{[\uu]}$ depends linearly on $[\uu]$;
\item
The map 
\begin{align*}
\h & \longrightarrow \R^{\Gamma_0} \\
[\uu] & \longmapsto \alpha_{[\uu]} 
\end{align*}
is injective (\cite{DrummGoldman2}).
\item
Suppose $\dim\vV=3$. 
If $\Gamma_{[\uu]}$ acts properly, then
$\vert\alpha_{[\uu]}(\gamma)\vert$ is the Lorentzian length of the
unique closed geodesic in $\eE/\Gamma_{[\uu]}$ corresponding to $\gamma$.
\end{itemize}
The significance of the {\em sign\/} of $\alpha_{[\uu]}$ is the following
result, due to Margulis~\cite{Margulis1,Margulis2}:

\begin{OSLemma*}
Suppose $\Gamma$ acts properly.
Either $\alpha_{[\uu]}(\gamma) > 0$ for all $\gamma\in\Gamma_0$ or
$\alpha_{[\uu]}(\gamma) < 0$ for all $\gamma\in\Gamma_0$.
\end{OSLemma*}

This paper provides a more conceptual proof of the Opposite Sign Lemma
by extending a normalized version of Margulis's invariant $\alpha$ 
to a continuous function $\Psi_\uu$ on a connected space $\gc$ of probability measures. 
If
\begin{equation*}
\alpha_{[\uu]}(\gamma_1)<0< \alpha_{[\uu]}(\gamma_2), 
\end{equation*}
then an element $\mu\in\gc$ ``between'' $\gamma_1$ and $\gamma_2$ exists,
satisfying $\Psi_\uu(\mu)=0$.

Associated to the linear group $\Gamma_0$ is a complete hyperbolic
surface $\Sigma = \Gamma_0\backslash \Ht$.  A {\em geodesic
current~\cite{Bonahon}\/} is a Borel probability measure on the unit
tangent bundle $\U$ of $\Sigma$ invariant under the geodesic
flow $\phi$.  We modify $\alpha_{[\uu]}$ by dividing it by the length
of the corresponding closed geodesic in the hyperbolic surface.  This
modification extends to a continuous function $\Psi_{[\uu]}$ 
on the compact space $\gc$ of {\em geodesic currents\/} on $\Sigma$.

Here is a brief description of the construction, which first
appeared in Labourie~\cite{Labourie}. For brevity we only describe
this for the case $\dim \vV=3$.
Corresponding to the affine deformation $\Gamma_{[\uu]}$ is 
a flat affine bundle $\E$ over $\U$, whose associated flat vector bundle
has a parallel Lorentzian structure $\B$.
Let $\phi$ be the vector field on $\U$ generating the geodesic flow.
For a sufficiently smooth section $s$ of $\E$, 
the covariant derivative $\nabla_\phi(s)$
is a section of the flat vector bundle $\V$ associated to $\E$.
Let $\nu$ denote the section of $\V$ which associates to a point
$x\in\U$ the unit-spacelike vector corresponding to the geodesic
associated to $x$. Let $\mu\in\gc$ be a geodesic current. Then
$\Psi_{[\uu]}(\mu)$ is defined by:
\begin{equation*}
\int_\U \B\big(\nabla_\phi(s),\nu\big) \, d\mu.
\end{equation*}

Nontrivial elements $\gamma\in\Gamma_0$ correspond to 
periodic orbits $c_\gamma$ for $\phi$.
The period of $c_\gamma$ equals the length
$\ell(\gamma)$ of the corresponding closed geodesic in $\Sigma$.
Denote the subset of $\gc$ consisting of measures supported
on periodic orbits by $\gcp$. 
Denote by $\mu_\gamma$ the geodesic current corresponding to 
$c_\gamma$.
Because $\alpha(\gamma^n) = \vert n\vert \alpha(\gamma)$ and 
$\ell(\gamma^n) = \vert n\vert \ell(\gamma)$, the ratio
$\alpha(\gamma)/\ell(\gamma)$ depends only on 
$\mu_\gamma$.

\begin{thm*}
Let $\Gamma_\uu$ denote an affine deformation of $\Gamma_0$. 
\begin{itemize}
\item
The function
\begin{align*}
\gcp & \longrightarrow \R \\
\mu_\gamma  & \longmapsto \frac{\alpha(\gamma)}{\ell(\gamma)}
\end{align*}
extends to a continuous function 
\begin{equation*}
\gc\xrightarrow{\Psi_{[\uu]}} \R. 
\end{equation*}
\item
$\Gamma_\uu$ acts properly if and only if 
$\Psi_{[\uu]}(\mu)\neq 0$ for all $\mu\in\gc$.
\end{itemize}
\end{thm*}
\noindent Since $\gc$ is connected, either $\Psi_{[\uu]}(\mu)> 0$ for all 
$\mu\in\gc$ or $\Psi_{[\uu]}(\mu)< 0$ for all
$\mu\in\gc$.
Compactness of $\gc$ implies that Margulis's invariant grows linearly with the
length of  $\gamma$:
\begin{corollary*}
For any $\uu\in\z$, there exists $C > 0$ such that
\begin{equation*}
\vert \alpha_{[\uu]}(\gamma)\vert \le C \ell(\gamma)
\end{equation*}
for all $\gamma\in \Gamma_\uu$.
\end{corollary*}

Margulis~\cite{Margulis1,Margulis2} originally used the invariant
$\alpha$ to detect properness: if $\Gamma$ acts properly, then all
the $\alpha(\gamma)$ are all positive, or all of them are negative.
(Since conjugation by a homothety scales the invariant: 
\begin{equation*}
\alpha_{\lambda \uu} = \lambda\alpha_{\uu}
\end{equation*}
uniform positivity and uniform negativity are equivalent.)
Goldman~\cite{Goldman, GoldmanMargulis}
conjectured this necessary condition is sufficient:
that is, properness is equivalent to the positivity (or negativity)
of Margulis's invariant. This has been proved 
when $\Sigma$ is a three-holed sphere; see
Jones~\cite{Jones} and Charette-Drumm-Goldman~\cite{CharetteDrummGoldman}.
In this case the cone corresponding to proper actions is defined 
by the inequalities corresponding to the three boundary components
$\gamma_1,\gamma_2,\gamma_3\subset \partial\Sigma$:
\begin{align*}
\big\{ [\uu]\in  & \h  \mid  \alpha_{[\uu]}(\gamma_i) > 0 \text{~for $i=1,2,3$}\big\} \\ 
& \bigcup\, \big\{ [\uu]\in \h \mid  \alpha_{[\uu]}(\gamma_i) < 0 \text{~for $i=1,2,3$}\big\}
\end{align*}
However, when $\Sigma$ is a one-holed torus,
Charette~\cite{Charette} showed that the positivity
of $\alpha_{[\uu]}$ requires infinitely many inequalities, 
suggesting the original conjecture is generally false.
Recently Goldman-Margulis-Minsky~\cite{GoldmanMargulisMinsky} have
disproved the original conjecture, using ideas inspired by Thurston's
unpublished manuscript~\cite{Thurston}.

This paper is organized as follows. \S 1 collects facts about convex cocompact hyperbolic surfaces and 
the recurrent part of the geodesic flow. \S 2 reviews constructions in affine geometry such as affine deformations,
and defines the Margulis invariant. \S 3 defines the flat affine bundles associated to an affine deformation. 
\S 4 introduces the horizontal lift of the geodesic flow to the flat affine bundle whose dynamics 
reflects the dynamics of the affine deformation of the Fuchsian group $\Gamma_0$. \S 5 extends the normalized
Margulis invariant to the function $\Psi_{[\uu]}$ on geodesic currents (first treated by Labourie~\cite{Labourie}).
\S 6 shows that for every {\em non-proper\/} affine deformation $[\uu]$ there is a geodesic current $\mu\in\gc$ for
which $\Psi_{[\uu]}(\mu) = 0$.
\S 7 shows that, conversely, for every {\em proper\/} affine deformation $[\uu]$, 
$\Psi_{[\uu]}(\mu) \neq 0$ for every $\mu\in\gc$. The Opposite Sign Lemma follows as a corollary.

\section*{Acknowledgements}
We thank Herbert Abels, Francis Bonahon, Virginie Charette, Todd
Drumm, David Fried, Yair Minsky, and Jonathan Rosenberg for helpful
conversations. We are grateful to the hospitality of the Centre
International de Recherches Math\'ematiques in Luminy where this paper
was initiated. We are also grateful to Gena Naskov
% g-noskov@yahoo.com
for several corrections, as well as the anonymous referee for helpful
comments.

\section*{Notation and Terminology}
Denote the tangent bundle of a smooth manifold $M$ by $TM$.
For a given Riemannian metric on $M$, {\em the unit tangent bundle\/}
$UM$ consists of unit vectors in $TM$.
For any subset $S\subset M$, denote the induced bundle over $S$ by $US$.
Denote the (real) hyperbolic plane by $\Ht$ and
its boundary by $\partial\Ht = \so$. 
For distinct $a,b\in\so$, denote the geodesic having $a,b$
as ideal endpoints by $\overleftrightarrow{ab}$.

Let $\V$ be a vector bundle over a manifold $M$. Denote by $\Az$ 
the vector space of sections of $\V$. 
Let $\nabla$ be a connection on $\V$. If $\ss\in\Az$ is a section and
$\xi$ is a vector field on $M$, then
\begin{equation*}
\nabla_\xi(\ss) \in \Az. 
\end{equation*}
denotes covariant derivative of $\ss$ with respect to $\xi$.
If $f$ is a function, then $\xi f$ denotes the directional derivative
of $f$ with respect to $\xi$.

An affine bundle $\E$ over $M$ is a fiber bundle 
over $M$ whose fiber is an affine space $\Ee$ with structure group
the group $\Aff$ of affine automorphisms of $\Ee$.
The  {\em vector bundle $\V$ associated to $\E$\/} is the vector bundle
whose fiber $\V_x$ over $x\in M$ is the vector space associated to the fiber
$\E_x$.  That is, $\V_x$ is the vector space consisting of
all translations $\E_x\longrightarrow \E_x$.

Denote the space of continuous sections of a bundle $\E$ by $\Gamma(\E)$. 
If $\V$ is a vector bundle, then $\Gamma(\V)$ is a vector space.
If $\E$ is an affine bundle with underlying vector bundle $\V$, 
then $\Gamma(\E)$ is an affine space with underlying vector space $\Gamma(\V)$.
If $s_1,s_2\in\Ss$ are two sections of $\E$, 
the {\em difference\/} $s_1-s_2$ is the section of $\V$, 
whose value at $x\in M$ is the translation $s_1(x) - s_2(x)$ of $\E_x$. 

Denote the convex set of Borel probability measures on a topological
space $X$ by $\Pr(X)$.  Denote the cohomology class of a cocycle $z$
by $[z]$.  A {\em flow\/} is an action of the additive group $\R$ of
real numbers.  The transformations in the flow defined by a vector
field $\xi$ are denoted $\xi_t$, for $t\in\R$.

\section{Hyperbolic geometry}\label{sec:hypgeo}

Let $\Rto$ be the 3-dimensional real vector space with inner product
\begin{equation*}
\bB (\vv,\ww) := \vv_1 \ww_1 + \vv_2 \ww_2 - \vv_3 \ww_3 
\end{equation*}
and $\Go  = \Oto$ the identity component of its group of isometries.
$\Go $ consists of linear isometries of $\Rto$ which preserve both an
orientation of $\Rto$ and a {\em time-orientation \/} (or {\em future\/})
of $\Rto$, that is, a connected component of the {\em open light-cone\/}
\begin{equation*}
\{ \vv\in\Rto \mid \bB(\vv,\vv) < 0 \}. 
\end{equation*}
Then 
\begin{equation*}
\Go  = \Oto \cong \PSL \cong \Isom^0 (\Ht)  
\end{equation*}
where $\Ht$ denotes the real hyperbolic plane.

\subsection{The hyperboloid model}
We define $\Ht$ as follows. We work in the irreducible representation
$\Rto$ (isomorphic to the adjoint representation of $\Go $).
Let $\bB$ denote the invariant bilinear form (isomorphic to the Killing
form). The two-sheeted hyperboloid
\begin{equation*}
\{ \vv\in \Rto \mid \bB(\vv,\vv) = -1\} 
\end{equation*}
has two connected components. If $\pp$ lies in this hyperboloid, then
its connected component equals
\begin{equation}\label{eq:component}
\{ \vv\in \Rto \mid \bB(\vv,\vv) = -1, \bB(\vv,\pp) <0 \}. 
\end{equation}
For clarity, we fix a timelike vector $\pp$, for example
\begin{equation*}
\pp =  \bmatrix 0 \\ 0 \\ 1 \endbmatrix \in \Rto, 
\end{equation*}
and define $\Ht$ as the connected component \eqref{eq:component}.

The Lorentzian metric defined by $\bB$ restricts to a Riemannian
metric of constant curvature $-1$ on $\Ht$. The identity component $\Go $
of $\Oto$ is the group $\Isomo$ of orientation-preserving isometries of $\Ht$.
The stabilizer of $\pp$ 
\begin{equation*}
\K := \mathsf{PO}(2).  
\end{equation*}
is the {\em maximal compact subgroup\/} of $\Go $.
Evaluation at $\pp$
\begin{align*}
\Go /\K & \longrightarrow \Rto \\
g \K & \longmapsto g \pp
\end{align*}
identifies $\Ht$ with the homogeneous space $\Go /\K$.

The action of $\Go $ canonically lifts to a {\em simply
transitive\/} (left-) action on the unit tangent bundle $U\Ht$.
The unit spacelike vector
\begin{equation}\label{eq:dotv}
\dot{\pp} := \bmatrix 0 \\ 1 \\ 0 \endbmatrix \in \Rto 
\end{equation}
defines a unit tangent vector to $\Ht$ at $\pp$, which we denote:
\begin{equation*}
(\pp,\dot{\pp}) \in U\Ht.
\end{equation*}
Evaluation at $(\pp,\dot{\pp})$
\begin{align}
\Go  & \stackrel{\Ee}\longrightarrow U\Ht \notag \\ 
g_0 & \longmapsto g_0 (\pp,\dot{\pp}) \label{eq:evaluation}
\end{align}
$\Go$-equivariantly identifies $\Go$ with the unit tangent bundle $U\Ht$.
Here $\Go$ acts on itself by left-multiplication.
The action on $U\Ht$ is the action induced by isometries of $\Ht$.
Under this identification, $\K$ corresponds 
to the fiber of $U\Ht$ above $\pp$.

\subsection{The geodesic flow}\label{sec:geodesicflow}

The {\em Cartan subgroup\/} 
\begin{equation*}
\A := \mathsf{PO}(1,1).  
\end{equation*}
is the one-parameter subgroup comprising
\begin{equation*}
a(t) = \bmatrix 1 & 0 & 0 \\
0 & \cosh(t) & \sinh(t) \\ 0 & \sinh(t) & \cosh(t) \endbmatrix
\end{equation*}
for $t\in\R$. Under this identification, 
\begin{equation*}
a(t) \pp \in \Ht
\end{equation*}
describes the geodesic through $\pp$ tangent to $\dot{\pp}$.

Right-multiplication by $a(-t)$ on $\Go $ identifies with the geodesic
flow $\tphi_t$ on $U\Ht$. First note that $\tphi_t$ on $U\Ht$ commutes
with the action of $\Go $ on $U\Ht$, which corresponds to the action on $\Go $
on itself by left-multiplications.  The $\R$-action on $\Go $ corresponding
to $\tphi_t$ commutes with left-multiplications on $\Go $.
Thus this $\R$-action on $\Go $ must be
right-multiplication by a one-parameter subgroup.  At the basepoint
$(\pp,\dot{\pp})$, the geodesic flow corresponds to 
right-multiplication by $\A$. Hence this one-parameter subgroup must be
$\A$. Therefore right-multiplication by $\A$ on
$\Go $ induces the geodesic flow $\tphi_t$ on $U\Ht$.

Denote the vector field corresponding to the geodesic flow by $\tphi$:
\begin{equation*}
\tphi(x) := \frac{d}{dt}\bigg|_{t=0}\tphi_t(x).
\end{equation*}

\subsection{The convex core}
Since $\Gamma_0$ is a Schottky group, the complete hyperbolic
surface is {\em convex cocompact.\/} That is, there exists a 
compact geodesically convex subsurface $\core$ such that the inclusion
\begin{equation*}
\core\subset \Sigma 
\end{equation*}
is a homotopy equivalence.
$\core$ is called the {\em convex core\/} of $\Sigma$, 
and is bounded by closed geodesics
$\partial_i(\Sigma)$ for $i=1,\dots,k$. %%% reference 
Equivalently the discrete subgroup $\Gamma_0$ is finitely generated
and contains only hyperbolic elements.

Another criterion for convex cocompactness involves the {\em
ends\/} $e_i(\Sigma)$ of $\Sigma$.  Each component $e_i(\Sigma)$ of
the closure of $\Sigma \setminus \core$ is diffeomorphic to a product
\begin{equation*}
e_i(\Sigma) \xrightarrow{\approx} \partial_i(\Sigma)\times [0,\infty).
\end{equation*}
The corresponding end of $\U$ is diffeomorphic to the product
\begin{equation*}
e_i(\U) \xrightarrow{\approx} \partial_i(\Sigma)\times S^1\times [0,\infty).
\end{equation*}
The corresponding Cartesian projections 
\begin{equation*}
e_i(\Sigma) \longrightarrow \partial_i(\Sigma)
\end{equation*}
and the identity map on $\core$ extend to a deformation
retraction 
\begin{equation}\label{eq:retraction} %%%%
\Sigma\xrightarrow{\Pi_\core}\core.  
\end{equation}
We use the following standard compactness results
(see Kapovich~\cite{Kapovich}, \S 4.17 (pp.\ 98--102), 
Canary-Epstein-Green~\cite{CEG},
or Marden~\cite{Marden} for discussion and proof):
%%% Katok-Hasselblatt, Brin-Stuck, Beford-Keane-Series???
\begin{lemma}\label{lem:NRcompact}
Let $\Sigma=\Ht/\Gamma_0$ be the hyperbolic surface where $\Gamma_0$ is
convex cocompact and let $R\ge 0$.
Then
\begin{equation}\label{eq:Rnbhd}
\kr := \{ y\in\Sigma \mid d(y,\core)\le R \} 
\end{equation}
is a compact geodesically convex subsurface, upon which the restriction of
$\Pi_\core$ is a homotopy equivalence.  
\end{lemma}

%\begin{proof}
%See Ratcliffe~\cite{Ratcliffe}.
%\end{proof}

\subsection{The recurrent subset}
Let $\Ur\subset \U$ denote the 
union of all recurrent orbits of $\phi$.

\begin{lemma}
Let $\Sigma$ be as above, $\U$ its unit tangent bundle and
$\phi_t$ the geodesic flow on $\U$.
Then:
\begin{itemize}
\item $\Ur\subset \U$ is a compact $\phi$-invariant
subset.
\item Every $\phi$-invariant Borel probability measure on $\U$ is
suppoted in $\Ur\subset \U$.
\item Tthe space of $\phi$-invariant Borel probability measures on $\U$ 
is a convex compact space $\gc$ with respect to the weak $\star$-topology.
\end{itemize}
\end{lemma}

\begin{proof}
Clearly $\Ur$ is invariant. 

The subbundle $U\core$ comprising 
unit vectors over points in $\core$ is compact.
Any geodesic which exits $\core$ into an end $e_i(\Sigma)$ 
remains in the same end $e_i(\Sigma)$.
Therefore every recurrent geodesic ray eventually lies in
$\core$. 

On the unit tangent bundle $U\Ht$, every $\phi$ orbit tends to
a unique ideal point on $\so =\partial\Ht$; let 
\begin{equation*}
U\Ht \xrightarrow{\eta} \so
\end{equation*}
denote the corresponding map. Then the $\phi$-orbit of $\tilde x\in U\Ht$
defines a  geodesic which is recurrent in the forward direction if and only
if $\eta(\tilde x)\in \Lambda$, where $\Lambda\subset\so$ 
is the {\em limit set\/}
of $\Gamma_0$. The geodesic is recurrent in the backward direction
if $\eta(-\tilde x)\in \Lambda$. Then every recurrent orbit of the geodesic
flow lies in the compact set consisting of vectors $x\in U\core$
with 
\begin{equation*}
\eta(\tilde x),\eta(-\tilde x)\in\Lambda. 
\end{equation*}
Thus $\Ur$ is compact.

By the Poincar\'e recurrence theorem~\cite{HasselblattKatok}, 
every $\phi$-invariant probability measure on $\U$ is supported in
$\Ur$. 
Since $\Ur$ is compact, 
the space of Borel probability measures
$\Pr(\Ur)$
on $\Ur$
is a compact convex set with respect to the weak topology.
The closed subset of $\phi$-invariant probability measures on $\Ur$ 
is also a compact convex set. Since
every $\phi$-invariant probability measure on $\U$ is supported
on $\Ur$, the assertion follows.
\end{proof}
The following well known fact will be needed later.

\begin{lemma} $\Ur$ is connected. \label{lem:connected}
\end{lemma}
\begin{proof} 
Let $\Lambda\subset\so$ be the limit set as above.
For every $\lambda\in\Lambda$, the orbit $\Gamma_0\lambda$
is dense in $\Lambda$. (See for example Marden~\cite{Marden},
Lemma 2.4.1.) Denote
\begin{equation*}
\Lambda^* :=  \{ (\lambda_1,\lambda_2)\in \Lambda\times\Lambda \mid  \lambda_1 \neq \lambda_2 \}.
\end{equation*}
The preimage $\tUr$ of $\Ur\subset\U$ under the
covering space $U\Ht\xrightarrow{\Pi} \U$ is the union
\begin{equation*}
\tUr \;=\; \bigcup_{(\lambda_1,\lambda_2)\in\Lambda^*} \overleftrightarrow{\lambda_1\lambda_2}
\end{equation*}
Choose a hyperbolic element $\gamma\in\Gamma_0$ and let $\lambda_\pm\in\Lambda$ denote its fixed points.
The image 
\begin{equation*}
c_\gamma \;:=\;  \Pi\big(\overleftrightarrow{\lambda_-\lambda+}\big)
\end{equation*}
is the unique closed geodesic in $\Ur$ corresponding to $\gamma$.

We claim that the image $U_+$ under $\Pi$ of the subset
\begin{equation*}
\bigcup_{\lambda\in\Lambda} \overleftrightarrow{\lambda_-\lambda} \;\subset\;\Ur
\end{equation*}
is connected. Suppose that $W_1,W_2\,\subset\,\U\,$ are disjoint open subsets such
that $U_+\subset W_1 \cup W_2$. 
Since $c_\gamma$ is connected, either $W_1$ or $W_2$ contains $c_\gamma$;
suppose that $W_1\supset c_\gamma$. We claim that $U_+\subset W_1$.
Otherwise $\overleftrightarrow{\lambda_-\lambda}$ meets $W_2$ for
some $\lambda\in\Lambda$.
Connectedness of  $\overleftrightarrow{\lambda_-\lambda}$ implies
$\overleftrightarrow{\lambda_-\lambda}\subset W_2$. Now
\begin{equation*}
\lim_{n\to\infty}\gamma^n\lambda = \lambda_+, 
\end{equation*}
which implies
\begin{equation*}
\lim_{n\to\infty}\gamma^n\big(\overleftrightarrow{\lambda_-\lambda}\big)
\;=\;
\lim_{n\to\infty}
\big(\overleftrightarrow{\lambda_-(\gamma^n\lambda)}\big)
\;=\;
\overleftrightarrow{\lambda_-\lambda_+}.
\end{equation*}
Thus for $n>>1$, the geodesic $\overleftrightarrow{\lambda_\lambda}$
lies in $\Pi^{-1}(W_1)$, a contradiction.

Thus $U_+\subset W_1$, implying that $U_+$ is connected.
Density of $\Gamma_0\lambda$ in $\Lambda$ implies 
$U_+$ is dense in $\Ur$. Therefore $\Ur$ is connected.
\end{proof}

%%%%%%%%%%%%%%%%%%%%%%%%%%
\section{Affine geometry }

This section collects general properties on affine spaces, affine
transformations affine deformations of linear group actions.
We are primarily interested in  affine deformations of linear
actions factoring through the irreducible $2r+1$-dimensional real
representation $\vV_r$ of 
\begin{equation*}
\Go \cong\PSL,  
\end{equation*}
where $r$ is a positive integer. 

\subsection{Affine spaces and their automorphisms}
Let $\vV$ be a real vector space. Denote its group of linear automorphisms 
by $\GLV$.

An {\em affine space \/} $\eE$ (modelled on $\vV$) 
is a space equipped with a simply transitive 
action of $\vV$. Call $\vV$ the {\em vector space underlying 
$\eE$,\/} and refer to its elements as {\em translations.\/}
Translation $\tau_\vv$ by a vector $\vv\in\vV$ 
is denoted by {\em addition:\/}
\begin{align*}
\eE & \xrightarrow{~\tau_{\vv}~} \eE \\
x & \longmapsto x + \vv. 
\end{align*}
Let $x,y\in \eE$.
Denote the unique vector $\vv\in\vV$ such that $\tau_\vv(x)=y$ by 
{\em subtraction:\/}
\begin{equation*}
\vv = y - x.  
\end{equation*}

Let $\eE$ be an affine space with associated vector space $\vV$. 
Choosing an arbitrary point $O\in\eE$ (the {\em origin\/}) identifies
$\eE$ with $\vV$ via the map
\begin{align*}
\vV &\longrightarrow \eE \\
\vv &\longmapsto O + \vv.
\end{align*}
An {\em affine automorphism\/}
of $\eE$ is the composition of a linear mapping 
(using the above identification
of $\eE$ with $\vV$) and a translation, that is,
\begin{align*}
\eE & \xrightarrow{g} \eE \\
O + \vv & \longmapsto O + \L(g)(\vv) + \uu(g)
\end{align*}
which we simply write as
\begin{equation*}
\vv \longmapsto \L(g)(\vv) + \uu(g). 
\end{equation*}
The affine automorphisms of $\eE$ form a group $\AffE$,  and 
$(\L,\uu)$ defines an isomorphism of $\AffE$ with 
the semidirect product %%$\vV\rtimes\GLV$.
$\GLV\ltimes\vV$.
The linear mapping $\L(g)\in\GLV$ 
is the {\em linear part\/} of the affine transformation $g$, 
and 
\begin{equation*}
\AffE \xrightarrow{\L} \GLV 
\end{equation*}
is a homomorphism.
The vector $\uu(g)\in \vV$ is the {\em translational part\/} of $g$.
The mapping
\begin{equation*}
\AffE \xrightarrow{\uu} \vV 
\end{equation*}
satisfies a {\em cocycle identity:\/}
\begin{equation}\label{eq:cocycle}
\uu(\gamma_1\gamma_2) = \uu(\gamma_1) + \L(\gamma_1) \uu(\gamma_2)
\end{equation}
for $\gamma_1,\gamma_2\in\AffE$.

\subsection{Affine deformations of linear actions}

Let $\Gamma_0\subset \GLV$ be a group of linear automorphisms of a
vector space $\vV$. Denote the corresponding $\Gamma_0$-module by
$\vV$ as well.

An {\em affine deformation \/} of $\Gamma_0$ is a representation 
\begin{equation*}
\Gamma_0 \stackrel{\rho}\longrightarrow \AffE 
\end{equation*}
such that $\L\circ\rho$ is the inclusion $\Gamma_0\hookrightarrow\GLV$.
We confuse $\rho$ with its image $\Gamma :=\rho(\Gamma_0)$, to which we
also refer as an {\em affine deformation\/} of $\Gamma_0$.
Note that $\rho$ embeds $\Gamma_0$ as the subgroup $\Gamma$
of $\GLV$.
In terms of the semidirect product decomposition
\begin{equation*}
\AffE \cong \vV \rtimes \GLV 
\end{equation*}
an affine deformation is the graph $\rho = \rho_\uu$
(with image denoted $\Gamma = \Gamma_\uu$) of a {\em cocycle\/}
\begin{equation*}
\Gamma_0 \stackrel{\uu}\longrightarrow V
\end{equation*}
that is, a map satisfying the cocycle identity \eqref{eq:cocycle}.
Write
\begin{equation*}
\gamma = \rho(\gamma_0) =  (\uu(\gamma_0),\gamma_0)\in V \rtimes \Gamma_0 
\end{equation*}
for the corresponding affine transformation:
\begin{equation*}
 x \stackrel{\gamma}\longmapsto \gamma_0 x + \uu(\gamma_0).
\end{equation*}

Cocycles form a vector space $\z$. 
Cocycles $\uu_1,\uu_2\in\z$ are {\em cohomologous\/}
if their difference $\uu_1-\uu_2$ is a {\em coboundary,\/} a cocycle
of the form
\begin{align*}
\Gamma_0 & \xrightarrow{\,\delta \vv_0\,} \vV \\
\gamma & \longmapsto  \vv_0 - \gamma \vv_0
\end{align*}
where $\vv_0\in\vV$.
Cohomology classes of cocycles form a vector space $\h$.  Affine
deformations $\rho_{\uu_1},\rho_{\uu_2}$ are conjugate by translation by
$\vv_0$ if and only if 
\begin{equation*}
\uu_1 - \uu_2 = \delta \vv_0. 
\end{equation*}
Thus $\h$ parametrizes translational conjugacy classes of affine deformations 
of $\Gamma_0\subset\GLV$

An important affine deformation of $\Gamma_0$ is the {\em trivial affine
deformation:\/} When $\uu=0$, the affine deformation 
$\Gamma_\uu$ equals $\Gamma_0$ itself. 

\subsection{Margulis's invariant of affine deformations}\label{sec:marginv}
Consider the case that $\Go =\PSL$ and $\L$ is an irreducible 
representation of $\Go $. For every positive integer $r$, let
$\L_r$ denote the irreducible representation of $\Go $ 
on the $2r$-symmetric power $\vV_r$ of the standard representation of 
$\SL$ on $\R^2$. The dimension of $\vV_r$ equals $2r+1$.
The central element $-\Id\in\SL$ acts by $(-1)^{2r}=1$, so this
representation of $\SL)$ defines a representation of 
\begin{equation*}
\PSL=\SL/\{\pm\Id\}. 
\end{equation*}
The representation $\Rto$ introduced in \S\ref{sec:hypgeo}
is $\vV_1$, the case when $r=1$.
Furthermore the $\Go $-invariant nondegenerate 
skew-symmetric bilinear form on $\R^2$ induces a nondegenerate 
{\em symmetric\/} bilinear form $\bB$ on $\vV_r$, which we normalize in
the following paragraph.

An element $\gamma\in \Go $ is {\em hyperbolic\/} if it corresponds to
an element $\tilde\gamma$ of $\SL$ with distinct real eigenvalues.
Necessarily these eigenvalues are reciprocals $\lambda,\lambda^{-1}$,
which we can uniquely specify by requiring $\vert\lambda\vert <1$.
Furthermore we choose eigenvectors $\vv_+,\vv_-\in\R^2$ 
such that:
\begin{itemize}
\item $\tilde\gamma(\vv_+) = \lambda\vv_+$;
\item $\tilde\gamma(\vv_-) = \lambda^{-1}\vv_-$;
\item The ordered basis $\{\vv_-,\vv_+\}$ is positively oriented.
\end{itemize}
Then the action $\L_r$ has eigenvalues $\lambda^{2j}$, for 
\begin{equation*}
j = -r, 1 -r,\dots -1,0,1,\dots,r-1, r,  
\end{equation*}
where the symmetric product
\begin{equation*}
\vv_-^{r-j}\vv_+^{r+j}\in \vV_r 
\end{equation*}
is an eigenvector with eigenvalue $\lambda^{2j}$.
In particular $\gamma$ fixes the vector 
\begin{equation*}
\xo(\gamma) := c\vv_-^r\vv_+^r, 
\end{equation*}
where the scalar $c$ is chosen so that
\begin{equation*}
\bB( \xo(\gamma),\xo(\gamma)) = 1.
\end{equation*}
Call $\xo(\gamma)$ the {\em neutral vector\/} of $\gamma$. 

The subspaces
\begin{align*}
\Vm(\gamma) &:= \,\sum_{j=1}^r \,\R \big(\vv_-^{r+j}\vv_+^{r-j}\big) \\
\Vp(\gamma) &:= \,\sum_{j=1}^r \, \R \big(\vv_-^{r-j}\vv_+^{r+j}\big)
\end{align*}
are $\gamma$-invariant and $\vV$
enjoys a $\gamma$-invariant 
$\bB$-orthogonal direct sum decomposition
\begin{equation*}
\vV \;=\; 
\Vm(\gamma) \,\oplus\, \R\big(\xo(\gamma)\big)  \,\oplus\, \Vp(\gamma).
\end{equation*}
\noindent
For any norm $\Vert\, \Vert$ on $\vV$, there exists $C,k>0$ such that
\begin{align}\label{eq:hyp}
\Vert \gamma^n (\vv) \Vert \le C e^{- k n} \Vert\vv\Vert &\text{for $\vv\in\Vp(\gamma)$} \notag \\
\Vert \gamma^{-n} (\vv) \Vert \le C e^{- k n} \Vert\vv\Vert &\text{for $\vv\in\Vm(\gamma)$}.
\end{align}
Furthermore
\begin{equation*}
\xo(\gamma^{n}) = \vert n\vert \xo(\gamma) 
\end{equation*}
if $n\in\Z, n\neq 0$, and
\begin{equation*}
\Vpm(\gamma^n) = \begin{cases} 
\Vpm(\gamma) &\text{~if $n>0$} \\
\Vmp(\gamma) &\text{~if $n<0$} \end{cases}.
\end{equation*}
For example, consider the hyperbolic one-parameter subgroup $\A$ 
comprising $a(t)$ where
\begin{equation*}
a(t): \begin{cases} \vv_+ & \longmapsto e^{t/2} \vv_+ \\
\vv_- & \longmapsto e^{-t/2} \vv_- \end{cases}\
\end{equation*}
In that case the action on $\vV_1$ corresponds to the one-parameter
group of isometries of $\Rto$ defined by
\begin{equation}\label{eq:alpha}
a(t) := \bmatrix 
\cosh(t) & 0 & \sinh(t) \\ 0 & 1 & 0 \\
\sinh(t) & 0 & \cosh(t) \endbmatrix
\end{equation}
with neutral vector
\begin{equation}\label{eq:neutral}
\xo\big(a(t)\big) = \bmatrix 0 \\ 1 \\0 \endbmatrix  
\end{equation}
when $t\neq 0$.

Next suppose that $g\in\AffE$ is an affine transformation whose linear
part $\gamma = \L(\gamma)$ is hyperbolic. Then there exists a unique
affine line $l_g\subset E$ which is $g$-invariant. The line $l_g$ is
parallel to the neutral vector $\xo(\gamma)$. The restriction of $g$
to $l_g$ is a translation by the vector
\begin{equation*}
\bB\big(\,g x - x\,,\, \xo(\gamma)\,\big)\;\; \xo(\gamma)
\end{equation*}
where $\xo(\gamma)$ is chosen so that $\bB(\xo(\gamma),\xo(\gamma))=1$.

Suppose that $\Gamma_0\subset \Go $ be a {\em Schottky group,\/} that is, a
nonabelian discrete subgroup containing only hyperbolic elements. Such
a discrete subgroup is a free group of rank at least two. The adjoint
representation of $\Go $ defines an isomorphism of $\Go $ with the
identity component $\Oto$ of the orthogonal group of the 3-dimensional
Lorentzian vector space $\Rto$. 
% Schottky group of 2nd kind? 

Let $\uu\in\z$ be a cocycle defining an affine deformation $\rho_\uu$
of $\Gamma_0$. In \cite{Margulis1,Margulis2}, 
Margulis constructed the invariant
\begin{align*}
\Gamma_0 & \xrightarrow{\alpha_\uu} \R \\
\gamma & \longmapsto  \bB(\uu(\gamma), \xo(\gamma)).
\end{align*}
This well-defined class function on $\Gamma_0$ satisfies
\begin{equation*}
\alpha_\uu(\gamma^n)  = \vert n\vert \alpha_\uu(\gamma) 
\end{equation*}
%% (when $r$ is odd)
%% should we include this???
%%
and depends only the cohomology class $[\uu]\in\h$.
Furthermore $\alpha_\uu(\gamma)=0$ if and only if $\rho_\uu(\gamma)$ fixes
a point in $E$. Two affine deformations of a given $\Gamma_0$ are conjugate
if and only if they have the same Margulis invariant 
(Drumm-Goldman~\cite{DrummGoldman2}). An affine deformation $\gamma_\uu$ of 
$\Gamma_0$ is {\em radiant\/} if it satisfies any of the following equivalent
conditions:
\begin{itemize}
\item $\Gamma_\uu$ fixes a point;
\item $\Gamma_\uu$ is conjugate to $\Gamma_0$;
\item The cohomology class $[\uu]\in\h$ is zero;
\item The Margulis invariant $\alpha_\uu$ is identically zero.
\end{itemize}
For further discussion see  \cite{Charette, CharetteDrumm,
DrummGoldman, Goldman, GoldmanHirsch, GoldmanMargulis,Labourie}.)

The centralizer of $\Gamma_0$ in the general linear group $\GLthr$ consists
of {\em homotheties\/} 
\begin{equation*}
\vv \stackrel{h_\lambda}\longmapsto  \lambda \vv
\end{equation*}
where $\lambda\neq 0$. The homothety $h_\lambda$ conjugates an affine
deformation $(\rho_0,\uu)$ to $(\rho_0,\lambda\uu)$. Thus conjugacy
classes of non-radiant affine deformations are parametrized by the
projective space $\Ph$. Our main result is that conjugacy classes
of proper actions comprise a convex domain in $\Ph$.

\section{Flat bundles associated to affine deformations}
We define two fiber bundles over $\U$.
The first bundle $\V$ is a {\em vector bundle\/} associated to the original
group $\Gamma_0$. The second bundle $\E$ is an {\em affine bundle\/}
associated to the affine deformation $\Gamma$. The vector bundle
underlying $\E$ is $\V$. The vector bundle $\V$ has a flat linear connection
and the affine bundle $\E$ has a flat affine connection, each denoted
$\nabla$. (For the theory of connections on affine bundles, see 
Kobayashi-Nomizu~\cite{KobayashiNomizu}.)
%% The vector field $\phi$ tangent to the flow $\Phi$ defined above 
%% I don't think $\Phi$ $
We define a vector field $\Phi$ on $\E$ which is uniquely determined
by the following properties:
\begin{itemize}
\item $\Phi$ is $\nabla$-horizontal ; %%
\item $\Phi$ covers the vector field $\phi$ defining the geodesic flow on $\U$.
\end{itemize}
We derive an direct alternate description of the corresponding flow $\Phi_t$
in terms of the Lie group $\Go$ and its semidirect products $\vV\rtimes\Go$.

%%%% this part should be completely rewritten!
\subsection{Semidirect products and homogeneous affine bundles}
Consider a Lie group $\Go $, a vector space $\vV$, and a linear
representation $\Go  \xrightarrow{\L}\GLV$. 
Let $\G$ be the corresponding semidirect product $\vV\rtimes \Go $.
Multiplication in $\G$ is defined by
\begin{equation}\label{eq:semidirect}
(\vv_1,g_1) \,\cdot\, (\vv_2, g_2)  \;:=\; (\vv_1 + g_1 \vv_2, g_1g_2). %%%%
\end{equation}
Projection 
\begin{align*}
\G & \xrightarrow{\Pi} \Go   \\
(\vv,g) &\longmapsto g
\end{align*}
defines a trivial bundle with fiber
$\vV$ over $\Go $. It is equivariant with respect to the action of
$\G$ on the total space by left-multiplication and the action of $\G$ on the
base obtained from left-multiplication on $\Go $ and the homomorphism $\L$. 
Since $\L$ is a homomorphism, \eqref{eq:semidirect} implies equivariance
of $\Pi$. 

When $r =1$, that is, $\vV = \Rto$, then $\G$ is the {\em tangent bundle\/}
$T\G_0$, with its natural Lie group structure. Compare Goldman-Margulis~\cite{GoldmanMargulis}
and \cite{Goldman}.

Since the fiber of $\Pi$ equals the vector space$\vV$, the fibration
$\G \xrightarrow{\Pi} \Go$ has the structure of (trivial) {\em affine
bundle\/} over $\Go$.  Furthermore this structure is $\G$-invariant:
For \eqref{eq:semidirect} implies that the action of $\gamma_1 =
(\vv_1,g_1)$ on the total space covers the action of $g_1$ on the
base. On the fibers the action is affine with linear part
$g_1=\L(\gamma_1)$ and translational part $\vv_1=\uu(\gamma_1)$.

Denote the total space of this $\G$-homogeneous affine bundle over $\Go$
by $\tE$.

Consider $\Pi$ also as a (trivial) {\em vector bundle.\/}
By \eqref{eq:semidirect}, this structure is $\Go$-invariant.
Via $\L$, this $\Go$-homogeneous vector bundle becomes
a $\G$-homogeneous vector bundle $\tV\longrightarrow\Go$.

{\em This  $\G$-homogeneous vector bundle underlies $\tE$:\/}
%%%% fix this!
Let
$(\vv_2'-\vv_2,1)$ be the translation taking   %%%
$(\vv_2,g_2)$ to $(\vv_2',g_2)$.
Then \eqref{eq:semidirect} implies that
$\gamma_1 = (\vv_1,g_1)$ acts on 
$(\vv_2'-\vv_2,g_2)$ by $\L$:
\begin{equation*}
\big( (\vv_1 + g_1(\vv_2')) - (\vv_1 + g_1(\vv_2)),g_1g_2\big) \;=\;
\big( \L(g_1)(\vv_2'-\vv_2),g_1g_2\big).
\end{equation*}

Of course both $\tE$ and $\tV$ identify with $\G$, 
but each has different actions of the discrete group $\Gamma \cong \Gamma_0$,
imparting the different structures of a {\em flat affine bundle\/} 
and a {\em flat vector bundle\/} to the respective quotients.

In our examples, $\L$ preserves a a bilinear form $\bB$ on $\vV$.
The $\Go $-invariant bilinear form $V\times V\xrightarrow{\bB}\R$ defines
a bilinear pairing $\tV\times\tV\xrightarrow{\B}\R$ of vector bundles.

\subsection{Homogeneous connections}
The $\G$-homogeneous affine bundle $\tE$
and the $\G$-homogeneous vector bundle $\tV$ admit flat connections 
(each denoted $\tilde\nabla$) as follows. To specify $\tilde\nabla$,
it suffices to define the covariant derivative of a section $\ts$
over a smooth path $g(t)$ in the base. For either $\tE=G$ or $\tV=G$, 
a section is determined by a smooth path $\vv(t)\in\vV$:
\begin{align*}
\R &\xrightarrow{~\ts~} \vV \rtimes \Go  = \G  \\
t &\longmapsto  (\vv(t),g(t)).
\end{align*}
Define
\begin{equation*}
\Ddt \ts(t) := \frac{d}{dt} \vv(t) \in \vV. 
\end{equation*}
If now $\ts$ is a section of $\tE$ or $\tV$, and $X$ is a tangent
vector field, define
\begin{equation*}
\tnabla_X(\ts) = \Ddt \ts(g(t)) 
\end{equation*}
where $g(t)$ is any smooth path with 
% \begin{equation*}
$
g'(t)=X\big(g(t)\big). 
$
% \end{equation*}
The resulting covariant differentiation operators define connections
on $\V$ and $\E$ which are invariant under the respective $\G$-actions.

\subsection{Flatness}
For each $\vv\in\vV$, 
\begin{align*}
\Go  &\xrightarrow{~\ts_\vv~} \vV\rtimes \Go  = \G \\
g &\longmapsto (\vv,g)
\end{align*}
defines a section whose image is the coset $\{\vv\} \times \Go  \subset \G$.
Clearly these sections are parallel with respect to $\tnabla$.
Since the sections $\ts_\vv$ foliate $\G$, the connections are flat. 

If $\L$ preserves a bilinear form $\bB$ on $\vV$, 
the bilinear pairing $\B$ on $\tV$ is parallel with respect to $\tnabla$.

\section{Sections and subbundles}
Now we describe the sections and subbundles of the homogeneous bundles
over $\Go \cong\PSL$ associated to the irreducible
$(2r+1)$-di\-men\-sional representation $\vV_r$.

\subsection{The flow on the affine bundle }
Right-multiplication by $a(-t)$ on $\G$ defines a flow $\tPhi_t$ on $\tE$.
Since $\A\subset \Go $, 
this flow covers the flow $\phi_t$ on $\Go $ 
defined by right-multiplication by $a(-t)$ on $\Go $,
where $a(-t)$ is defined by \eqref{eq:alpha}.
That is, the diagram
\begin{equation*}
\begin{CD} 
\tE @>{\tPhi_t}>>\tE  \\
@V{\tPi}VV @VV{\tPi}V \\
U\Ht @>>{\tphi_t}> U\Ht
\end{CD} 
\end{equation*}
commutes. The vector field $\tPhi$ on $\tE$ generating $\tPhi_t$ covers
the vector field $\tphi$ generating $\tphi_t$.
%%% changed \txi to \tphi

Furthermore %%%% $\ts_v(g a(-t))=\ts(v)a(-t)$ 
\begin{align*}
\ts_\vv\big(g a(-t)\big) & = 
\big(\vv, g a(-t)\big)  \\
& = (\vv,g) \big(0, a(-t)\big)  \\
& = \ts_\vv(g) a(-t)
\end{align*}
implies
%%% $\tphi_t\circ\ts_v = \ts_v\circ\tPhi_t$, 
\begin{equation*}
\ts_\vv\circ\tphi_t = \tPhi_t\circ\ts_\vv,    %%%% fixed v's
\end{equation*}
%%% correction
whence $\tPhi$ is the {\em $\tnabla$-horizontal
lift\/} of $\tphi$.

The flow $\tPhi_t$ commutes with the action of $\G$. Thus {\em
$\tPhi_t$ is a flow on the flat $\G$-homogeneous affine bundle % $\tV$
$\tE$
covering $\phi_t$.}

Right-multiplication by $a(t)$ on $\G$ also defines a flow on the flat
$\G$-homogeneous vector bundle $\tV$ covering $\tphi_t$. 
Identifying $\tV$ as the vector bundle underlying $\tE$, the $\R$-action is
just the linearization
$D\tPhi_t$ of the action $\tPhi_t$:
\begin{equation*}\begin{CD}
\tV @>{D\tPhi_t}>> \tV \\
@VVV @VVV \\
U\Ht @>>{\tphi_t}> U\Ht 
\end{CD}\end{equation*}

\subsection{The neutral section}
The $\G$-action and the flow $D\tPhi_t$ on $\tV$ preserve
a section $\tnu\in\tV$ defined as follows.
The one-parameter subgroup $\A$ fixes the neutral vector 
$\vv_0\in\vV$
defined in \eqref{eq:neutral}.
Let $\tnu$ denote the section of $\tV$ defined by:
\begin{align*}
U\Ht \approx \Go  &\xrightarrow{\tilde{\nu}} \vV\rtimes \Go  \approx \tV \\
g &\longmapsto (g \vv_0, g). %%% removed \L
\end{align*}
In terms of the group operation on $\G$, this section is given by {\em
right-multiplication\/} by $\vv_0\in V\subset \G$ acting on $g\in
\Go \subset \G$. Since $\L$ is a homomorphism,
\begin{equation*}
\tnu( h g) = h \tnu(g),
\end{equation*}
so $\tnu$ defines a $\G$-invariant section of $\tV$.

Although $\tnu$ is not parallel in every direction, 
it is parallel along the flow $\tphi_t$:
\begin{lemma}\label{lem:xinu}
$\tnabla_{\tphi}(\tnu) = 0$.
\end{lemma}
\begin{proof}
Let $g\in \Go $. Then
\begin{equation*}
\tnabla_{\tphi}(\tnu) (g)  
= \frac{D}{dt}\bigg|_{t=0} \tnu(\tphi_t(g)) =
\frac{D}{dt}\bigg|_{t=0} g a(-t) \vv_0 = 0
\end{equation*}
since $a(-t)\vv_0 = \vv_0$ is constant.
\end{proof}
\noindent
The section $\tnu$ is the diffuse analogue of the neutral eigenvector
$\xo(\gamma)$ of a hyperbolic element $\gamma\in\Oto$ discussed in
\S\ref{sec:marginv}. Another approach to the
flat connection $\nabla$ and the neutral section $\nu$ is given
in Labourie~\cite{Labourie}.

\subsection{Stable and unstable subbundles}
Let $\vV^0\subset\vV$ denote the line of vectors fixed by $\A$,
that is, the line spanned by $\vv_0$.
The eigenvalues  of $a(t)$ acting on $\vV$ are the $2r+1$ distinct positive
real numbers
\begin{equation*}
e^{rt}, e^{(r-1)t},\dots, 1, \dots,e^{(1-r)t},e^{-rt}   %% 
\end{equation*}
and the eigenspace decomposition of $V$ is invariant under $a(t)$.
Let $\Vp$ denote the sum of eigenspaces for eigenvalues $>1$ and
$\Vm$ denote the sum of eigenspaces for eigenvalues $<1$. The corresponding
decomposition
\begin{equation}\label{eq:decomposeV}
\vV = \Vm \oplus \vV^0 \oplus \Vp 
\end{equation}
is $a(t)$-invariant and defines a (left)%%
$\G$-invariant decomposition of the
vector bundle $\tV$ into subspaces which are invariant under $D\tPhi_t$.

%%%%%%%%%%%%%%%%%%%%%%%%%%%%%%%%%%%%%%%%% big section removed here
\subsection{Bundles over $\U$}\label{sec:bundles}
Let $\Gamma\subset \G$ be an affine deformation of a discrete subgroup 
$\Gamma_0\subset \Go $. Since $\Gamma$ is a discrete subgroup of $\G$, the
quotient $\E := \tE/\Gamma$ is an affine bundle over 
$\U=U\Ht/\Gamma_0$ and inherits a flat connection $\nabla$ from
the flat connection $\tnabla$ on $\tE$.
Furthermore the flow $\tPhi_t$ on $\tE$ descends to a flow $\Phi_t$
on $\E$ which is the horizontal lift of the flow $\phi_t$ on $\U$.

The vector bundle $\V$ underlying $\E$ is the quotient
\begin{equation*}
\V := \tV/\Gamma = \tV/\Gamma_0 
\end{equation*}
and inherits a flat linear connection
$\nabla$ from the flat linear connection $\tnabla$ on $\tV$.
The flow $D\tPhi_t$ on $\tV$ covering $\tphi_t$, the neutral section
$\tnu$, and the stable-unstable splitting
\eqref{eq:decomposeV}
all descend to a
flow $D\Phi_t$, a section $\nu$ and a splitting
\begin{equation}\label{eq:splitting}
\V = \V^- \oplus \V^0 \oplus \V^+ 
\end{equation}
of the flat vector bundle $\V$ over $\U$.
%  In particular
% \begin{equation*}
% \big(\Pi_\core\big)^* \V_\core \cong \V
% \end{equation*}
% since $\Pi_\core$ is a homotopy equivalence.
%%%% do we actually need this? clarify 
%%%% if we do, refer to retraction

There exists a Euclidean metric $\gg$ on $\V$ with the following properties:
\begin{itemize}
\item The neutral section $\nu$ is bounded with respect to $\gg$;
\item The bilinear form $\B$ is bounded with respect to $\gg$;
\item The flat linear connection $\nabla$ is bounded with respect to $\gg$;
\item  {\em Hyperbolicity:\/}
The flow $D\Phi_t$ exponentially expands the subbundle $\Vp$
and exponentially contracts $\Vm$. 
Explicitly,
there exist 
constants $C,k>0$ exist so that 
\begin{equation}\label{eq:hyperbolicityplus}
\Vert D(\Phi_t) \vv \Vert \le C e^{ kt} \Vert\vv \Vert 
\text{\quad as $t\longrightarrow-\infty$},
\end{equation}
for $\vv\in\V^+$, and
\begin{equation}\label{eq:hyperbolicityminus}
\Vert D(\Phi_t)\vv \Vert \le C e^{ -kt} \Vert\vv \Vert 
\text{\quad as $t\longrightarrow+\infty$}
\end{equation}
for $\vv\in\V^-$.
%% which follow immediately from \eqref{eq:hyp}.
\end{itemize}
To construct such a metric, first choose 
a Euclidean inner product on the vector space $\vV$.
Then extend it to a left $\Go$-invariant 
Euclidean metric on the vector bundle 
\begin{equation*}
\tV\approx \G = \vV \rtimes \Go.  
\end{equation*}
Hyperbolicity follows from the description of the
adjoint action of $\A$ as in \eqref{eq:hyp}.
%% which is invariant under left-multiplications in $\Go$.
%% This Euclidean metric satisfies 
%% the following {\em hyperbolicity condition:\/}
%%
With this Euclidean metric, the space $\Gamma(\V)$ of continuous
sections of $\V$ is a Banach space.

% \subsection
\section{Proper $\Gamma$-actions and proper $\R$-actions}

Let $X$ be a locally compact Hausdorff space with homeomorphism
group $\homeo$. Suppose that  $H$ is a closed subgroup 
of $\homeo$ with respect to the compact-open topology.
Recall that $H$ acts {\em properly\/} if the mapping
\begin{align*}
H \times X & \longrightarrow X \times X \\
(g, x) & \longmapsto (gx,x) 
\end{align*}
is a proper mapping. That is, for every pair of compact subsets
$K_1,K_2\subset X$, the set
\begin{equation*}
\{ g\in H \mid gK_1\cap K_2\neq\emptyset  \}
\end{equation*}
is compact. The usual notion of a {\em properly discontinuous\/} action
is a proper action where $H$ is given the discrete topology.
In this case, the quotient $X/H$ is a Hausdorff space. If $H$ acts
freely, then the quotient mapping $X \longrightarrow X/H$ is a covering
space. For background on proper actions, see
Bourbaki~\cite{Bourbaki}, Koszul~\cite{Koszul} and Palais~\cite{Palais}.

The question of properness of the $\Gamma$-action is equivalent
to that of properness of an action of $\R$.
\begin{proposition}\label{prop:prop}
An affine deformation $\Gamma$ acts properly on $E$ if and only if
$\A$ acts properly by right-multiplication on $\G/\Gamma$.
\end{proposition}

The following lemma is a key tool in our argument. 
(A related statement is proved in Benoist~\cite{Benoist}, Lemma~3.1.1. 
For a proof in a different context see Rieffel~\cite{Rieffel1,Rieffel2}.)

\begin{lemma}\label{lem:properness}
Let $X$ be a locally compact Hausdorff space. Let
$A$ and $B$ be commuting groups of homeomorphisms of $X$, each of
which acts properly on $X$. 
Then the following conditions are equivalent:
\begin{enumerate}
\item $A$ acts properly on $X/B$;\label{it:A}
\item $B$ acts properly on $X/A$;\label{it:B}
\item $A\times B$ acts properly on $X$.\label{it:AB}
\end{enumerate}
\end{lemma}
\begin{proof}
We prove (\ref{it:AB}) $\Longrightarrow$ (\ref{it:A}).
Suppose  $A\times B$ acts properly on $X$ but 
$A$ does not act properly on $X/B$. Then 
a $B$-invariant subset $H\subset X$ exists  such that:
\begin{itemize}  
\item $H/B\subset X/B$ is compact; 
\item  The subset
\begin{equation*}
A_{H/B} := \{ a\in A \mid a(H/B) \cap (H/B) \neq \emptyset\} 
\end{equation*}
is not compact.
\end{itemize}
% By Koszul~\cite{Koszul} or Palais~\cite{Palais}, 
We claim a compact subset 
$K\subset X$ exists satisfying $B\cdot K = H$.

Denote the quotient mapping by
\begin{equation*}
X \xrightarrow{\Pi_B} X/B. 
\end{equation*}
For each $x\in\Pi_B^{-1}(H/B)$, choose a
precompact open neighborhood $U(x)\subset X$. The images $\Pi_B(U(x))$ define
an open covering of $H/B$. Choose a finite subset $x_1,\dots, x_l$ such that
the open sets $\Pi_B(U(x_i))$ cover $H/B$.
Then the union \begin{equation*}
K' := \bigcup_{i=1}^l \overline{U(x_i)} 
\end{equation*}
is a compact subset of $X$ such that $\Pi_B(K')\supset H/B$.
Taking $K = K'\cap \Pi_B^{-1}(H/B)$ proves the claim.

Since $A$ acts properly on $X$, the subset
\begin{equation*}
(A\times B)_{K} := \{ (a,b)\in A\times B \mid aK \cap bK \neq \emptyset\} 
\end{equation*}
is compact. However $A_{H/B}$ is the image of the compact set
$(A\times B)_{K}$ under Cartesian projection $A\times B\longrightarrow A$
and is compact, a contradiction.

We prove that (\ref{it:A}) $\Longrightarrow$ (\ref{it:AB}).
Suppose that $A$ acts properly on $X/B$ and  
$K\subset X$ is compact. Cartesian projection 
$A\times B\longrightarrow A$ maps
\begin{equation}\label{eq:cartesianprojection}
(A\times B)_K \longrightarrow  A_{(B\cdot K)/B}. 
\end{equation}
$A$ acts properly on $X/B$ implies $A_{(B\cdot K)/B}$ is
a compact subset of $A$.  Because \eqref{eq:cartesianprojection} is a 
proper map, $(A\times B)_K$ is compact, as desired.

Thus (\ref{it:AB}) $\Longleftrightarrow$  (\ref{it:A}) .
The proof (\ref{it:AB}) $\Longleftrightarrow$ (\ref{it:B}) is similar.

\end{proof}
\begin{proof}[Proof of Proposition~\ref{prop:prop}]
Apply Lemma~\ref{lem:properness} 
with $X=\G$ and $A$ to be the action of %%
$A_0 \cong\R$ by right-multiplication and $B$
to be the action of $\Gamma$ by left-multiplication. The lemma implies
that $\Gamma$ acts properly on $\eE = \G/\Go $ if and only if $\Go $ acts
properly on $\Gamma\backslash \G$.

Apply the Cartan decomposition 
\begin{equation*}
\Go  = \K \A \K. 
\end{equation*}
Since $\K$ is compact,  the action of $\Go $
on $E$ is proper if and only if its restriction to $\A$ is proper.
\end{proof}

\section{Labourie's diffusion of Margulis's invariant}
In \cite{Labourie}, Labourie
defined a function 
\begin{equation*}
\U\xrightarrow{F_\uu}\R,  
\end{equation*}
corresponding to the invariant $\alpha_\uu$ defined by
Margulis~\cite{Margulis1,Margulis2}.  Margulis's
invariant $\alpha = \alpha_\uu$ is an $\R$-valued class function on
$\Gamma_0$ whose value on $\gamma\in \Gamma_0$ equals
\begin{equation*}
B\big(\,\rho_\uu(\gamma) O - O,  \xo(\gamma)\,\big)
\end{equation*}
where $O$ is the origin and $\xo(\gamma)\in\vV$ is the neutral vector of
$\gamma$ (see \S\ref{sec:marginv} and the references listed there).

Now the origin $O\in E$ will be replaced by a section $s$ of $\E$, 
the vector $\xo(\gamma)\in\vV$ will be replaced by the neutral section 
$\nu$ of $\V$, 
and the linear action of $\Gamma_0$ on $\vV$ will be replaced by the geodesic flow 
$\phi_t$ on $\U$. 

Let $s$ be a 
$C^\infty$ section of $\E$.
%% continuous section of $\E$ %%
%% which is continuously differentiable along $\phi$. 
In particular $s$ is $C^1$ along $\phi$. That is, the function
\begin{align*}
\R &\longrightarrow \E \\
t &\longmapsto s\big(\phi_t(p)\big) 
\end{align*}
is $C^1$ for all $p\in \U$. 
%%
%  I think we need smooth sections, since we differentiate and subtract 
% in section 8.
%%%%
%%%%
Its covariant derivative with respect to $\phi$
is a smooth section $\nabla_{\phi}(s)\in\Az$. 
Pairing with $\nu\in\Az$ via 
\begin{equation*}
\V\times\V\xrightarrow{\B}\R 
\end{equation*}
produces a  continuous function 
$\U\xrightarrow{\F}\R$ defined by: 
\begin{equation*} 
\F := \B( {%% \tilde
\nabla}_\phi(s),\nu).
\end{equation*}

\subsection{The invariant is continuous}

Let $\Ss$ denote the space of continuous  sections  $s$ of $\E$ over 
$\Ur$  which are differentiable along $\phi$  and the covariant 
derivative  $\nabla_\phi(s)$ is continuous.
If $s\in\Ss$, then $\F$ is continuous.

%Henceforth choose a Riemannian metric $\gg$ on $\V$ such that
%$\bB$ and $\nu$ are bounded with respect to $\gg$.
%Let $\Ss$ denote the space of 
%
%sections $s$ of $\E$ satisfying the
%following properties:
%\begin{itemize}
%\item $s$ is continuous;
%\item $s$ is $C^1$ along $\phi$;
%\item $\nabla_\phi(s)$ is bounded with respect to $\gg$.
%\end{itemize} 
%In particular $\F$ is bounded and %% measurable.
%%
For each probability measure $\mu\in\Pr(\U)$, 
\begin{equation*}
\int_{\U} \F\, d\mu  
\end{equation*}
is a well-defined real number and the function
\begin{align*}
\Pr(\U) & \longrightarrow \R \\
\mu & \longmapsto \int_{\U} \F\, d\mu  
\end{align*}
is continuous in the weak $\star$-topology
on $\Pr(\U)$. 
Furthermore its restriction to the 
subspace $\gc \subset \Pr(\U)$ of
$\Phi$-invariant measures is also
continuous in the weak $\star$-topology.

\subsection{The invariant is independent of the section}

\begin{lemma}\label{lem:independent}
Let $s_1,s_2\in\Ss$ and $\mu$ be a $\phi$-invariant Borel probability measure 
on $\U$.
Then 
\begin{equation*}
\int_{\U} F_{\uu,s_1}\, d\mu  = \int_{\U} F_{\uu,s_2}\, d\mu.  
\end{equation*}
\end{lemma}
\begin{proof}
The difference $\ss = s_1 - s_2 $ of the sections $s_1,s_2$ of the
affine bundle $\E$ is a section of the vector bundle $\V$ and
\begin{equation*}
\nabla_\phi(s_1) -\nabla_\phi(s_2) =
\nabla_\phi(\ss).
\end{equation*}
Therefore
\begin{align*}
\int_{\U} F_{\uu,s_1}\, d\mu  - \int_{\U} F_{\uu,s_2} \,d\mu 
&  = \int_{\U} \B( {\nabla}_\phi(\ss),\nu) \,d\mu \\
& = \int_{\U} \phi \B( \ss,\nu) \,d\mu  -
\int_{\U} \B( \ss,\nabla_\phi(\nu)) \,d\mu.
\end{align*}
The first term vanishes since $\mu$ is $\phi$-invariant:
\begin{align*}
\int_{\U} \phi \B( \ss,\nu) \,d\mu  & =  
\int_{\U} \frac{d}{dt} (\phi_t)^*\big(\B( \ss,\nu)\big) \,d\mu  \\ & =  
\frac{d}{dt} \int_{\U} (\phi_t)^*\big(\B( \ss,\nu)\big) \,d\mu  \\ & =  
\frac{d}{dt} \int_{\U} \B( \ss,\nu) d\big((\phi_t)_*\mu\big)  = 0.  
\end{align*}
The second term vanishes since $\nabla_\phi(\nu)=0$ (Lemma~\ref{lem:xinu}).
\end{proof}
\noindent Thus 
\begin{equation*}
\Psi_{[\uu]}(\mu) := \int_{\U} \F \,d\mu  
\end{equation*}
is a well-defined continuous function
$\gc \xrightarrow{\Psi_{[\uu]}} \R$ 
which is independent of the section $s$ used to define it.
%% whose continuity follows from the definition of the weak topology
%% on $\Pr(\U)$.

\subsection{Periodic orbits}
Suppose 
\begin{equation*}
t \longmapsto \phi_t(x_0)
\end{equation*}
defines a periodic orbit of the geodesic
flow $\phi$ with {\em period\/} $T > 0$.
That is, $T$ is the least positive number
such that $\phi_{T}(x_0) = x_0$.
Suppose that $\gamma$ is the corresponding element of 
$\Gamma_0 \cong \pi_1(\Sigma)$.
Then $T$ equals the length $\ell(\gamma)$ of the 
corresponding closed geodesic on $\Sigma$.

%In terms of the identification $\U \approx \Gamma_0\backslash\Go$,
%write:
%\begin{equation*}
%x_0 = \Gamma_0 g_0 \in \Gamma_0 \backslash \Go = \U
%\end{equation*}
%for some $g_0\in\Go$ and
%\begin{equation*}
%\phi_t(x_0) =  \Gamma_0  g_0a(-t). 
%\end{equation*}
%Then $\gamma$ is the unique element of $\Go$
%for which
%\begin{equation*}
%\gamma g_0 a(-T)  = g_0,
%\end{equation*}
%that is,
%\begin{equation*}
%\gamma = g_0 a(T) g_0^{-1}.
%\end{equation*}
%
%\begin{equation*} %%
%\Gamma_0  ga(t + T)
%\;=\;  \Gamma_0 \gamma g a(t) \;=\;  \Gamma_0 g a(t)
%\end{equation*}
%where $T = \ell(\gamma)$ is the period of the orbit.

%%%% useful notation, maybe move it between 4.1 and 4.2
%%%% if it will be used to discuss covariant differentiation
%%%% and parallel transport

% The following notation is useful. 
If $x_0\in\U$ and $T > 0$, let 
$\phi^{x_0}_{[0,T]}$ denote the map
\begin{align}\label{eq:orbitseg}
[0,T] & \,\xrightarrow{~~\phi^{x_0}_{[0,T]}~}\, \U \notag \\
t & \longmapsto \phi_t(x_0).
\end{align}
The period  $T$ of the periodic orbit corresponding to $\gamma$
equals the length $\ell(\gamma)$ of the closed geodesic in $\Sigma$
corresponding to $\gamma$.
Then
\begin{equation}\label{eq:periodic}
\mu_\gamma := \frac1{T} (\phi^{x_0}_{[0,T]})_*\big(\mu_{[0,T]}\big)
\end{equation}
defines the {\em geodesic current associated to the periodic 
orbit $\gamma$,\/}
where $\mu_{[0,T]}$ denotes Lebesgue measure on $[0,T]$.

\begin{proposition}[Labourie~\cite{Labourie}, Proposition~4.2]
\label{prop:labourie}
Let $\gamma\in\Gamma_0$ be hyperbolic and let $\mu_\gamma\in\gcp$ be
the corresponding geodesic current. Then
\begin{equation}\label{eq:closedorbit}
\alpha(\gamma) =  \ell(\gamma) \int_\U \,\F\, d\mu_\gamma.
\end{equation}
\end{proposition}
\begin{proof}

Let $x_0\in\U$ is a point on the periodic orbit, and consider the 
map \eqref{eq:orbitseg} and the geodesic current \eqref{eq:periodic}.

Choose a section $s\in\Ss$. 
Pull $s$ back by the covering space 
% \begin{equation*}
% \begin{CD} 
% \tE @>{\tPhi_t}>>\tE  \\
% @V{\tPi}VV @VV{\tPi}V \\
% U\Ht @>>{\tphi_t}> U\Ht
% \end{CD} 
% \end{equation*}
\begin{equation*}
U\Ht\approx \Go \xrightarrow{\Pi} \U  
\end{equation*}
to a section of $\tE = \Pi^*\E$.
Since $\tE\longrightarrow U\Ht$ is a trivial $\Ee$-bundle,
this section 
is the graph of a map
\begin{equation*}
\Go  \xrightarrow{\tvv} \eE 
\end{equation*}
satisfying
\begin{equation*}
\tvv\circ \gamma =  \rho(\gamma) \tvv.
\end{equation*}
Then 
\begin{equation}\label{eq:integaloverclosedorbit}
\int_\U \,\F\, d\mu_\gamma =
\frac1{T} \int_0^T \F(\phi_t(x_0)) dt,
\end{equation}
which we evaluate by passing to the covering space $U\Ht\approx \Go $.

Lift the basepoint $x_0$ to $\tx_0\in U\Ht$. 
Then $\tx_0$ corresponds via $\Ee$ 
to some $g_0\in \Go$, 
where $\Ee$ is defined in \eqref{eq:evaluation}. 
The path $\phi^{x_0}_{[0,T]}$ lifts 
to the map 
\begin{align*}
[0,T] & \xrightarrow{\tphi^{\tx_0}_{[0,T]}} U\Ht \\
t & \longmapsto \tphi_t(\tx_0) 
\longleftrightarrow  g_0 a(-t).
\end{align*}
The periodic orbit lifts to the trajectory
\begin{align*}
\R &\longrightarrow \Go  \\
t &\longmapsto g_0 a(-t).
\end{align*}
Since the periodic orbit corresponds to the deck transformation $\gamma$
(which acts by left-multiplication on $\Go $),
\begin{equation*}
\gamma g_0 = g_0 a(-T) 
\end{equation*}
which implies
\begin{equation}\label{eq:aNegativet}
\gamma = g_0 a(-T) g_0^{-1}.
\end{equation}
\noindent 
Evaluate $\nabla_\phi s$ and $\nu$ along the trajectory
$\tphi^{\tx_0}_{[0,T]}$
by lifting to the covering space and computing in $\G$:
% Lift $s$ to 
% \begin{align*}
% \Go  & \xrightarrow{\ts} V \rtimes \Go  = G \\
% g &\longmapsto (\tvv(g),g)
% \end{align*}
% where $\tvv:\Go \longrightarrow V$ satisfies
% \begin{equation*} %\label{eq:equivariance}
% \tvv(\gamma g) = \rho(\gamma)\tvv(g).
% \end{equation*}
\begin{equation}\label{eq:evaluateCDeriv}
\big(\nabla_\phi\ts\big)(\phi_t\tx_0) = \Ddt \tvv(g_0 a(-t)).
\end{equation}
\noindent
In semidirect product coordinates, the lift $\tnu$ is defined by the
map
\begin{align}
\Go  &\xrightarrow{\tnu} \vV \notag\\
g &\longmapsto \L(g) \vv_0. \label{eq:evaluateNeutral}
\end{align}
%\end{proof}
\begin{lemma}\label{lem:neutralPer}
For any $g\in\Go$ and $t\in \R$, 
\begin{equation*}
\tnu\big(g a(-t)\big)= \xo(\gamma). 
\end{equation*}
\end{lemma}
\begin{proof}
\begin{align*}
\tnu\big(g a(-t)\big) & = \L\big(g a(-t)\big)\; \vv_0 
\quad\qquad \text{~\big(by \eqref{eq:evaluateNeutral}\big)}
\\ & =
\L\big(g a(-t)\big)\; \xo\big(a(-T)\big) \\ & =
\xo\big( (g a(-t)) a(-T) (g a(-t))^{-1}\big) \\ & = 
\xo\big( (g a(-T) g^{-1}\big)  = \xo(\gamma) 
\quad\qquad\text{ ~(by \eqref{eq:aNegativet})}.
\end{align*}
\end{proof}

%\begin{proof}[Conclusion of proof of Proposition~\ref{prop:labourie}]
\noindent{\it Conclusion of proof of Proposition~\ref{prop:labourie}.\/}
% \newline
\noindent
Evaluate the integrand in \eqref{eq:integaloverclosedorbit}:
\begin{align*}
\F(\phi_t x_0) & = \B(\nabla_\phi s,\nu) (\phi_t x_0) \\ &  =
\B\big(\nabla_\phi\ts(g_0 a(-t)), \tnu(g_0 a(-t)\big) \\ & = 
\bB\bigg(\Ddt \tvv(g_0 a(-t)), \xo(\gamma)\bigg) 
\text{~\big(by Lemma~\ref{lem:neutralPer} and \eqref{eq:evaluateCDeriv}\big)}
\\ & = 
\ddt \bB\big(\tvv(g_0 a(-t)), \xo(\gamma)\big) 
\end{align*}
and
\begin{align*}
\int_0^T \F(\phi_t x_0) dt & = 
\bB\big(\tvv(g_0 a(-T)), \xo(\gamma)\big) - \bB\big(\tvv(g_0), \xo(\gamma)\big) \\ 
& = \bB\big( \tvv(\gamma g_0) - \tvv(g_0), \xo(\gamma)\big) \\ 
& = \bB\big( \rho(\gamma)\tvv(g_0) - \tvv(g_0), \xo(\gamma)\big)  \\
& = \alpha(\gamma)
\end{align*}
as claimed.
\end{proof}

\section{Nonproper deformations}
Now we prove that $0 \notin \Psi_{[\uu]}(\gc)$ implies 
$\Gamma_\uu$ acts properly.

\begin{proposition}\label{prop:properness}
Suppose that $\Gamma_\uu$ is a non-proper affine deformation. Then there exists a
geodesic current  $\mu\in\gc $(that is, a $\phi$-invariant Borel probability measure)
such that
$\Psi_{[\uu]}(\mu)=0$.
\end{proposition}
\begin{proof}
By Proposition~\ref{prop:prop}, we may assume that the flow $\Phi_t$ defines a
non-proper action on $\E$. Choose compact subsets
$K_1,K_2\subset \E$ for which the set of $t\in\R$ such that
\begin{equation*}
\Phi_t(K_1)\cap K_2 \neq \emptyset 
\end{equation*}
is noncompact. 
Choose an unbounded sequence $t_n\in\R$, and sequences
$P_n\in K_1$ and $Q_n\in K_2$ such that 
\begin{equation*}
\Phi_{t_n}P_n = Q_n. 
\end{equation*}
Passing to subsequences, assume that 
$t_n \nearrow +\infty$ and 
\begin{equation}\label{eq:PQlimit}
\lim_{n\to\infty} P_n = P_\infty, \qquad 
\lim_{n\to\infty} Q_n = Q_\infty.
\end{equation}
for some $P_\infty,Q_\infty\in \E$.
For $n=1,\dots,\infty$, the images 
\begin{equation*}
p_n=\Pi(P_n),\; q_n=\Pi(Q_n) 
\end{equation*}
are points in $\U$ such that
\begin{equation}\label{eq:pqlimit}
\lim_{n\to\infty}p_n  = p_\infty,\qquad  
\lim_{n\to\infty}q_n  = q_\infty 
\end{equation}
and $\phi_{t_n}(p_n)=q_n$.
Choose $R>0$ such that 
\begin{align*}
d(p_\infty,U\core)\, &<\; R, \\ 
d(q_\infty,U\core)\, &<\; R. 
\end{align*}
Passing to a subsequence assume that all
$p_n,q_n$ lie in the compact set $U\kr$ %%
where $\kr$ is the $R$-neighborhood of the
core, defined in \eqref{eq:Rnbhd}. 
Since $\kr$ is geodesically convex, the curves
\begin{equation*}
\{ \phi_t(p_n) \mid 0\le t \le t_n \} 
\end{equation*}
also lie in $U\kr$(Lemma~\ref{lem:NRcompact}).

Choose a section $s\in\Ss$.
We use the splitting \eqref{eq:splitting} of $\V$
and the section $s$ of $\E$ to decompose 
the points $P_1,\dots, P_\infty$ and 
$Q_1,\dots, Q_\infty$.
For $n=1,\dots,\infty$, write
\begin{align}\label{eq:decomposepq}
P_n & = s(p_n) + \big( \ppn^- + a_n \nu + \ppn^+ \big) \notag \\
Q_n & = s(q_n) + \big( \qqn^- + b_n \nu + \qqn^+ \big)
\end{align}
where $a_n, b_n\in \R$ , $\ppn^- , \qqn^- \in \V^-$, and 
$\ppn^+,\qqn^+ \in \V^+$.
Since 
\begin{equation*}
D\Phi_t(\nu)=\nu 
\end{equation*}
and $\Phi_{t_n}(P_n)=Q_n$, taking $\V^0$-components in
\eqref{eq:decomposepq} yields:
\begin{equation}\label{eq:intab}
\int_0^{t_n} (\F\circ \phi_t)(p_n)\, dt \;=\; b_n - a_n.
\end{equation}
Since $s$ is continuous, \eqref{eq:pqlimit} implies
$s(p_n)\to s(p_\infty)$ and $s(q_n)\to s(q_\infty)$. 
By \eqref{eq:PQlimit},
\begin{equation*}
\lim_{n\longrightarrow\infty} (b_n - a_n) \;=\; b_\infty - a_\infty.
\end{equation*}
Thus \eqref{eq:intab} implies
\begin{equation*}
\lim_{n\to\infty} 
\int_0^{t_n} (\F\circ\phi_t)(p_n)\, dt \;=\; 
b_\infty - a_\infty 
\end{equation*}
so
\begin{equation}\label{eq:convergence}
\lim_{n\to\infty}  \frac1{t_n}
\int_0^{t_n} (\F\circ\phi_t)(p_n)\, dt \;=\; 0
\end{equation}
(because $t_n\nearrow+\infty)$.

Now we construct the measure $\mu\in\gc$ such that
\begin{equation*}
\Psi_{[\uu]}(\mu)  \;=\; \int_{\U} \F\, d\mu  \;=\; 0. 
\end{equation*}
Define approximating probability measures 
\begin{equation*}
\mu_n\in\Pr(U\kr) 
\end{equation*}
by pushing forward Lebesgue measure on the orbit through $p_n$ and
dividing by $t_n$:
\begin{equation*}
\mu_n := \,\frac1{t_n}\, (\phi^{p_n}_{[0,t_n]})_* \mu_{[0,t_n]}
\end{equation*}
where $\phi^{p_n}_{[0,t_n]}$ is defined in \eqref{eq:orbitseg}.
Compactness of $U\kr$ 
guarantees a weakly convergent subsequence 
$\mu_n$ in $\Pr(U\kr)$. 
Let
\begin{equation*}
\mu_\infty := \lim_{n\longrightarrow\infty} \mu_n.
\end{equation*}
Thus by \eqref{eq:convergence}.
\begin{equation*}
\Psi_{[\uu]}(\mu_\infty) =  \lim_{n\to\infty}  \frac1{t_n}
\int_0^{t_n} (\F\circ\phi_t)(p_n) dt = 0.
\end{equation*}

Finally we show that $\mu_\infty$ is $\phi$-invariant. 
For any $f\in L^\infty(U\kr)$ and $\lambda\in\R$,
\begin{equation*}
\bigg\vert \int f d\mu_n - \int (f\circ \phi_\lambda) d\mu_n \bigg\vert <
\frac{2\lambda}{t_n} \Vert f\Vert_\infty 
\end{equation*}
for $n$ sufficiently large. Passing to the limit,
\begin{equation*}
\bigg\vert \int f \,d\mu_\infty \,-\, \int (f\circ \phi_\lambda) \,d\mu_\infty 
\bigg\vert  \;=\; 0 
\end{equation*}
as desired.
\end{proof}

\section{Proper deformations}
Now we prove that $\Psi_{[\uu]}(\mu)\neq 0$ for a proper deformation
$\Gamma_\uu$.  The proof uses a lemma ensuring that the section $s$ can
be chosen only to vary in the neutral direction $\nu$. The proof uses
the hyperbolicity of the geodesic flow.  The uniform positivity or
negativity of $\Psi_{[\uu]}$ implies Margulis's ``Opposite Sign
Lemma' (Margulis'~\cite{Margulis1,Margulis2}, Abels~\cite{Abels},
Drumm~\cite{Drumm3,Drumm5}) as a corollary.

\begin{proposition}\label{prop:nonproperness}
Suppose that $\Phi$ defines a proper action. Then 
\begin{equation*}
\Psi_{[\uu]}(\mu)\neq 0 
\end{equation*}
for all $\mu\in\gc$.
\end{proposition}

\begin{corollary}[Margulis~\cite{Margulis1,Margulis2}]
\label{cor:oppositesign} 
Suppose that $\gamma_1,\gamma_2\in\Gamma_0$ satisfy 
\begin{equation*}
\alpha(\gamma_1) < 0 < \alpha(\gamma_2). 
\end{equation*}
Then $\Gamma$ does not act properly.
\end{corollary}
\begin{proof}
Proposition~\ref{prop:labourie} implies
\begin{equation*}
\Psi_{[\uu]}(\mu_{\gamma_1}) < 0 < \Psi_{[\uu]}(\mu_{\gamma_2})
\end{equation*}
Convexity of $\gc$ implies a continuous path $\mu_t\in\gc$ exists, with 
$t\in [1,2]$, 
for which 
\begin{align*}
\mu_1 &= \mu_{\gamma_1}, \\ \mu_2 &= \mu_{\gamma_2}.  
\end{align*}
The function 
\begin{equation*}
\gc\xrightarrow{\Psi_{[\uu]}} \R
\end{equation*}
is continuous.
The intermediate value theorem implies $\Psi_{[\uu]}(\mu_t)=0$ for some 
$1 < t < 2$. Proposition~\ref{prop:nonproperness} 
implies that $\Gamma$ does not act properly.
\end{proof}

\subsection{Neutralizing sections}
A section $s\in\Ss$ is {\em neutralized\/} if and only if 
\begin{equation*}
\nabla_{\phi}(s) \in \V^0.
\end{equation*}
This will enable us to relate the properness of the flow 
$\{\Phi_t\}_{t\in\R}$ to the invariant
$\Psi_{[\uu]}(\mu)$. 
To construct neutralized sections, we need the following
technical facts. 

The flow $\{D\Phi_t\}_{t\in\R}$ on $\V$ and the flow
 $\{\Phi_t\}_{t\in\R}$ on $\E$ induce (by push-forward) one-parameter
groups of continuous bounded operators $(D\Phi_t)_*$
and $(\Phi_t)_*$ on
$\Gamma(\V)$ and $\Gamma(\E)$ respectively.
We begin with an elementary observation,
whose proof is immediate. (Compare \S\ref{sec:bundles}.)

\begin{lemma}\label{lem:deriv}
Let $\xi\in\Gamma(\V)$. 
Suppose that 
\begin{equation*}
t \longmapsto (D\Phi_t)_*(\xi)
\end{equation*}
is a path in the Banach space $\Gamma(\V)$ which is
differentiable at $t = 0$. 
Then $\xi\in\Ss$ and 
\begin{equation*}
\frac{d}{dt}\bigg\vert_{t=0}  (D\Phi_t)_*(\xi) = \nabla_\phi(\xi).
\end{equation*}
\end{lemma}

Here is our main lemma:

\begin{lemma}
A neutralized section exists.
\end{lemma}
\begin{proof}
By the hyperbolicity condition \eqref{eq:hyperbolicityplus},

%the operator norm $\Vert(D\Phi_t)_*\Vert_+$  of  the restriction
%of $(D\Phi_t)_*$ to $\Gamma(\V^+)$ satisfies
\begin{equation*}
\Vert (D\Phi_t)_*(\xi^+)\Vert_+ \le C e^{-k t} \Vert\xi^+\Vert
\end{equation*}
for any section $\xi^+\in\Gamma(\V^+)$.
Consequently the improper integral
\begin{equation*}
\zeta^+ := \int_0^\infty (D\Phi_t)_*(\xi^+) dt
\end{equation*}
converges. Moreover, Lemma~\ref{lem:deriv} implies that
$\zeta^+\in\Ss$ with  $$\nabla_\phi(\zeta^+) = - \xi^+.$$
Similarly \eqref{eq:hyperbolicityminus} implies that,
for any section $\xi^-\in\Gamma(\V^-)$, the improper integral
\begin{equation*}
\zeta^-:= \int_{-\infty}^0 (D\Phi_t)_*(\xi^-) dt
\end{equation*}
converges, and 
$\zeta^-\in\Ss$ with  $$\nabla_\phi(\zeta^-) =  \xi^-.$$
\noindent
Now let $s\in\Ss$.  Decompose $\nabla_\phi(s)$ by the splitting
\eqref{eq:splitting}: %%% into components:
\begin{equation*}
\nabla_{\phi}(s) = \nabla^-_{\phi}(s) + \nabla^0_{\phi}(s) + \nabla^+_{\phi}(s),
\end{equation*}
where $\nabla^\pm_{\phi}(s)\in\Gamma(\V^\pm)$. Apply the above
discussion to $\xi^- = \nabla^-_\phi(s)$ and $\xi^+ = \nabla^+_\phi(s)$.
Then
\begin{equation*}
s_0=s \;+\; 
\int_0^\infty  (D\Phi_t)_*\big(\nabla^-_\phi(s)\big)\,dt 
\; -\;  
\int_0^\infty  (D\Phi_{-t})_*\big(\nabla^+_\phi(s)\big)\,dt  %%% corrected this
\end{equation*}
lies in $\Ss$ with   $\nabla_\phi(s_0) = \nabla^0_\phi(s)$.
Hence $s_0$ is neutralized, as claimed.
\end{proof}

\subsection{Conclusion of the proof}
Let $\mu\in\gc$ so that 
\begin{equation*}
\Psi_{[\uu]}(\mu) = \int \F d\mu = 0.
\end{equation*}
Define, for $T> 0$ and $p\in \Ur$,
\begin{equation*}
g_T(p) :=   \int_0^T \F (\phi_tp) dt.
\end{equation*}
Since $\mu$ is $\phi$-invariant, by Fubini's Theorem,
\begin{equation*}
\int g_T d\mu = 0.
\end{equation*}
Therefore, for every $T>0$, since $\Ur$ is connected by Lemma ~\ref{lem:connected},  there exists
$p_T\in\Ur$ such that
\begin{equation*}
g_T(p_T) = 0.
\end{equation*}
We may assume that 
$s\in\Ss$ is neutralized.
Then
\begin{equation*}
\frac{d}{dt} (D\Phi_t)_*(s) = \nabla_\phi(s) = \F \nu,
\end{equation*}
and
\begin{equation}\label{eq:flow}
(\Phi_T s)(p) = s(\phi_T p) + \bigg( \int_0^T \F (\phi_tp) dt \bigg) \nu
\end{equation}
for any $T>0$. Thus
\begin{equation*}
(\Phi_T s)(p_T) = s(\phi_T p_T). 
\end{equation*}
Let $K$ be the compact set $s(\Ur)$. Then, for all $T>0$,
\begin{equation*}
\Phi_T(K) \cap K \neq \emptyset,
\end{equation*}
and $\{\Phi_t\}_{t\in\R}$ is not proper, as claimed.
This completes the proof.

\makeatletter \renewcommand{\@biblabel}[1]{\hfill#1.}\makeatother

\end{document}